\documentclass[12pt]{article}
\usepackage[english]{babel}
\usepackage{amsfonts}
\usepackage{amsmath}
\usepackage{amssymb}
\usepackage{mathrsfs}
\usepackage{latexsym}
\usepackage{multicol}
\usepackage{multirow}
\usepackage{fancybox} 
\usepackage{color}
\usepackage{hyperref}
\usepackage{graphicx} 
\usepackage{caption}
\usepackage{enumitem}
\usepackage{upgreek}
\usepackage[overload]{empheq}
\usepackage{float}
\date{}

\def\be{\begin{equation}}
\def\ee{\end{equation}}

\def\d'{``}

\newtheorem{thm}{Theorem}[section]

\newtheorem{propn}[thm]{Proposition}

\newtheorem{cor}[thm]{Corollary}

\def\be{\begin{equation}}
\def\ee{\end{equation}}
\def\bea{\begin{eqnarray}}
\def\eea{\end{eqnarray}}

\def\i'{\textrm{i}}

\def\d'{``}

\begin{document}

\title{Integral representations and zeros of the Lommel function and the hypergeometric $_1F_2$ function.} 
\author{Federico Zullo\thanks{DICATAM,  Universit\`a degli Studi di Brescia, Brescia, Italy \& INFN, \& INFN, Sezione di Milano-Bicocca, Milano, Italy}}

\maketitle 

\begin{abstract}
We give different integral representations of the Lommel function $s_{\mu,\nu}(z)$ involving trigonometric and hypergeometric $_2F_1$ functions. By using classical results of P\'olya, we give the distribution of the zeros of $s_{\mu,\nu}(z)$ for certain regions in the plane $(\mu,\nu)$. Further, thanks to a well known relation between the functions $s_{\mu,\nu}(z)$ and the hypergeometric  $ _1F_2$  function, we describe the distribution of the zeros of $_1F_2$ for specific values of its parameters.
\end{abstract}




\section{Introduction}\label{ra_sec1}
The Lommel function $s_{\mu,\nu}(z)$ is a particular solution of the inhomogeneous Bessel differential equation:
\begin{equation}\label{Lom}
z^2\frac{d^2y}{dz^2}+z\frac{dy}{dz}+(z^2-\nu^2)y=z^{\mu+1}.
\end{equation}
More precisely,  $s_{\mu,\nu}(z)$ is the solution of equation (\ref{Lom}) satisfying  $s_{\mu,\nu}(z)\sim \frac{z^{\mu+1}}{(\mu+1)^2-\nu^2}(1+O(z^2))$ around $z=0$ when $\mu\pm\nu$ is not an odd negative integer. The explicit Taylor series of $s_{\mu,\nu}(z)$ around $z=0$ is given by \cite{W}
\begin{equation}\label{sser}
s_{\mu,\nu}(z)=\frac{z^{\mu+1}}{(\mu+1)^2-\nu^2}\left(1-\frac{z^{2}}{(\mu+3)^2-\nu^2}+\frac{z^{4}}{\left((\mu+3)^2-\nu^2\right)\left((\mu+5)^2-\nu^2\right)}+\cdots\right)
\end{equation}
The previous expression can be equivalently written in terms of hypergeometric $_1F_2$ function as \cite{W}
\begin{equation}\label{iper12}
s_{\mu,\nu}(z)=\frac{z^{\mu+1}}{(\mu+1)^2-\nu^2}\,_1F_2\left(1;\frac{\mu-\nu+3}{2},\frac{\mu+\nu+3}{2};-\frac{z^2}{4}\right)
\end{equation}
Notice that the function $_1F_2$ in (\ref{iper12}), or $z^{-\mu-1}s_{\mu,\nu}(z)$, is an entire function of $z$ with order of growth equal to 1. 
The Lommel function has many applications in mathematical physics and applied sciences: for example, it is fundamental in the description of  piezoelectricity phenomena \cite{ASS}, optics \cite{BW}, radiative transfer processes in the atmosphere \cite{Cha}, mechanics \cite{Sitz}, elastic scattering theory \cite{Thom}. 

In this paper we will analyze the distribution of the zeros of the function $s_{\mu,\nu}(z)$. The identification of the values $(\mu,\nu)$ giving a non negative function $s_{\mu,\nu}(z)$ on the real axis and the description of the location of the zeros have been investigated by many authors. In \cite{Cooke} Cooke showed  that if $\nu\geq 0$ and $\nu\leq \mu-1$ then $s_{\mu,\nu}(z)>0$ for $z>0$. Equally, if $\nu\geq 1/2$ and $\nu\leq \mu$ $s_{\mu,\nu}(z)>0$ for $z>0$, unless $\mu=\nu=1/2$ when $s_{\mu,\nu}\geq 0$ for $z>0$.  Successively, Steinig \cite{Steinig} showed that $s_{\mu,\nu}(z)>0$ for $z>0$ if $|\nu|<\mu$ and $\mu \geq 1/2$ except when $|\nu|=\mu=1/2$, where one has $s_{\mu,\nu}(z)\geq 0$ for $z>0$. Also, Steinig established the following properties of $s_{\mu,\nu}(z)$: if $\mu <1/2$ or if $\mu=1/2$ and $|\nu|>1/2$ then  $s_{\mu,\nu}(z)$ has an infinite number of zeros for $z\in (0,\infty)$. An extension of the positivity results of Cooke and Steinig has been given more recently by Cho \& Chung \cite{Cho}, since they showed that $s_{\mu,\nu}>0$ also in the region $\mu>1/2$ and $\nu^2 \leq (\mu+1)^2-2$, whereas $s_{\mu,\nu}\geq 0$ for $\mu=|\nu|=1/2$, giving also interlacing properties among the zeros of $s_{\mu,\nu}$ and those of the Bessel function $J_{\nu}$ for certain values of the parameters $(\mu,\nu)$. In \cite{KL} Koumandos and Lamprecht give estimates about the location of the zeros of $s_{\mu,\nu}(z)$ for $\mu \in (-1/2,1/2)$ and $\nu=1/2$, whereas  in \cite{K} Koumandos, extends these to $\mu \in (-3/2,1/2)$, $\mu \neq -1/2$ and $\nu=1/2$. For more properties of the Lommel functios the reader can look for example at \cite{NIST}, Ch. 11.9 or \cite{W}, Ch. 10.7.

One of the tools that we will use is a result due to P\'olya \cite{P} about
the zeros of entire functions possessing suitable integral representations. The theorem is the following:
\begin{thm}\label{ThP}
[P\'olya, 1918 \cite{P}] Suppose that the function $f(t)$ is positive and not decreasing in $(0,1)$. Then the functions of $z$ defined by
\begin{equation}\label{P}
V(z)=\int_0^1\sin(zt)f(t)dt, \quad U(z)=\int_0^1\cos(zt)f(t)dt
\end{equation}
possesses only real zeros. Further,  if $f(t)$ grows steadily these zeros are simple and the intervals $(k\pi, (k+1)\pi)$, $k=1, 2, \ldots$ contain the positive zeros of $V(z)$, whereas the intervals $(\frac{(2k+1)\pi}{2}, \frac{(2k+3)\pi}{2})$, $k=0, 1, 2, \ldots$ contain the zeros of $U(z)$, each interval containing just one zero in both cases.
\end{thm}
By growing steadily it is meant that the function is not piecewise constant with a finite number of rational points of discontinuity in $(0,1)$.

\section{Some integral representations and comments}
In this section we will report some known facts about integral representations of the function $s_{\mu,\nu}(z)$. In particular, in 1936 Szymanski \cite{S} got an integral representation for the Lommel function in the case $\nu-\mu$ is an even positive integer. Actually,
Szymanski firstly considers equation (\ref{Lom}) for imaginary argument: 
\begin{equation}\label{LomIm}
z^2\frac{d^2T}{dz^2}+z\frac{dT}{dz}-(z^2+\nu^2)T=z^{\mu+1}.
\end{equation}
and noticed that the solution of (\ref{LomIm}) satisfying $T=t_{\mu,\nu}(z)\sim \frac{z^{\mu+1}}{(\mu+1)^2-\nu^2}(1+O(z))$ can be expressed as
\begin{equation}\label{S}
t_{\mu,\nu}=z^{\mu+1}\int_0^1\cosh(zt)W_{\mu,\nu}(t)dt,
\end{equation}
where $W_{\mu,\nu}$ is the solution of the following differential equation
\begin{equation}\label{W}
(1-t^2)\frac{d^2W_{\mu,\nu}}{dt^2}+(2\mu-1)t\frac{dW_{\mu,\nu}}{dt}+(\nu^2-\mu^2)W_{\mu,\nu}=0,
\end{equation}
with the following boundary conditions
\begin{equation}\label{Wbc}
W(1)=0, \quad W'(0)=-1.
\end{equation}
The solution of equation (\ref{W}) with the boundary conditions (\ref{Wbc}) is written by Szymanski in terms of the Gegenbauer functions $C_{n}^k$:
\begin{equation}\label{G}
W_{\mu,\nu}(t)=-(1-t^2)^{\mu+\frac{1}{2}}\frac{C_{\nu-\mu-1}^{\mu+1}(t)}{\left(C_{\nu-\mu-1}^{\mu+1}(t)\right)'_{t=0}}.
\end{equation}
where, for integers values of $\nu-\mu-1$ (odd values, for even positive integers the denominator vanishes), the functions $C_{n}^k$ are defined by the generating function $(1-2xt+x^2)^{-k}$, i.e. $C_n^k(t)$ is the coefficients of $x^n$ in the expansion of $(1-2xt+x^2)$ in powers of $x$:
\begin{equation}
(1-2xt+x^2)^{-k}=\sum_{n=0}C_n^k(t)x^n.
\end{equation} 
For the solution $s_{\mu,\nu}(z)$ of (\ref{Lom}) equally one finds, for $\nu-\mu$ even integers:
\begin{equation}\label{S}
s_{\mu,\nu}=z^{\mu+1}\int_0^1\cos(zt)W_{\mu,\nu}(t)dt,
\end{equation}
where $W_{\mu,\nu}(t)$ is again given by the expression (\ref{G}). 

In \cite{Z} another integral representation of the function $s_{0,v}(z)$ is given. In particular, by a direct check, it is possible to show that the function
\begin{equation}\label{int}
s_{0,\nu}(z)=\frac{1}{1+\cos(\pi\nu)}\int_0^\pi \sin(z\sin(t))\cos(\nu t)dt.
\end{equation}
solves equation (\ref{Lom}) when $\mu=0$ with  $s_{0,\nu}(z)\sim \frac{z}{1-\nu^2}(1+O(z))$. Equivalently, the function $s_{0,\nu}(z)$ can be also represented as 
\begin{equation}\label{int1}
s_{0,\nu}(z)=\frac{1}{1+\cos(\pi\nu)}\int_0^1 \sin(z t)\frac{\cos(\nu\arcsin(t))+\cos(\nu\pi-\nu\arcsin(t))}{\sqrt{1-t^2}}dt,
\end{equation}
where the range of $\arcsin$ is $(-\pi/2,\pi/2)$. It is possible to show that the function multiplying $\sin(zt)$ in (\ref{int1}) is positive and increasing for $|\nu|<1$ and $t\in (0,1)$. It is then possible to apply the Theorem (\ref{ThP}) to get the following

\begin{cor}\label{cor0}
The function $s_{0,\nu}(z)$, for $|\nu|<1$, possesses only real zeros. The zeros are simple and the intervals $(k\pi, (k+1)\pi)$, $k=0, 1, 2, \ldots$ contain the non-negative zeros of $s_{0,\nu}(z)$, each interval containing just one zero.
\end{cor}
In the next section we will extend the integral representation (\ref{int}) and the results described in Corollary (\ref{cor0}) to other values of $\mu$. Actually, the existence of a formula generalizing (\ref{int}), given by formula (\ref{intr}) below, has been pointed out to us after the writing of this work: it is a special case of an equality given by Prudnikov et al. \cite{Pru}. For the sake of completeness, we re-derive such formula in a manner that is closer to what Szymanski did in \cite{S} and to what is reported in formula (\ref{int1}). Some comments on the formula given in \cite{Pru} will be also given. The application to this formula to get a generalization of (\ref{cor0}) seems instead to be new to us.

\section{Other integral representations and the zeros}
In this section we will assume that $\mu\pm\nu$ is not an odd negative integer so that the series (\ref{sser}) is defined. By generalizing the case $\mu=0$ considered in \cite{Z}, we make the following ansatz for the function $s_{\mu,\nu}(z)$:
\begin{equation}\label{ans}
s_{\mu,\nu}(z)=z^{\mu}\int_0^1 \sin(zt)f_{\mu,\nu}(t)dt.
\end{equation}
We assume that the functions $f_{\mu,\nu}(t)$ are finite in $z=0$ whereas may be unbounded for $t\to 1$. By inserting (\ref{ans}) in (\ref{Lom}) we get that the functions $f_{\mu,\nu}(t)$ must be the solutions of the following differential equations
\begin{equation}\label{eqf}
(1-t^2)\frac{d^2f_{\mu,\nu}}{dt^2}+(2\mu-3)t\frac{df_{\mu,\nu}}{dt}+(\nu^2-(\mu-1)^2)f_{\mu,\nu}=0,
\end{equation}
with the boundary conditions:
\begin{equation}\label{bou}
f_{\mu,\nu}(0)=1,\quad \lim_{t\to 1^{-}}(1-t)f_{\mu,\nu}(t)=0, \quad \lim_{t\to 1^-}\left((1-t^2)\frac{df_{\mu,\nu}}{dt}+(2\mu-1)tf_{\mu,\nu}\right)=0.
\end{equation}
The general solution of equation (\ref{eqf}) can be represented in terms of hypergeometric functions as:
\begin{equation}
A(1-t)^{\mu-\frac{1}{2}}\, _2F_1\left(\frac{1}{2}+\nu,\frac{1}{2}-\nu;\mu+\frac{1}{2};\frac{1-t}{2}\right)+B\, _2F_1\left(1-\nu-\mu,1+\nu-\mu;\frac{3}{2}-\mu;\frac{1-t}{2}\right)
\end{equation}
where $A$ and $B$ are two arbitrary constants. The condition $ \lim_{t\to 1^{-}}(1-t)f_{\mu,\nu}(t)=0$ gives
\begin{equation}
\mu+\frac{1}{2}>0,
\end{equation}
whereas the series around $t=1$ of the function $(1-t^2)\frac{df_{\mu,\nu}}{dt}+(2\mu-1)tf_{\mu,\nu}$ is given by
\begin{equation}
(1-t^2)\frac{df_{\mu,\nu}}{dt}+(2\mu-1)tf_{\mu,\nu} \sim \left(A(1-t)^{\mu+\frac{1}{2}}\frac{\nu^2-\mu^2}{2\mu+1}+B(2\mu+1)\right)\left(1+O(1-t)\right),
\end{equation}
and results in $B=0$. Finally, from the boundary condition (\ref{bou}) $f_{\mu,\nu}(0)=1$ we get 
\begin{equation}
A\, _2F_1\left(\frac{1}{2}+\nu,\frac{1}{2}-\nu;\mu+\frac{1}{2};\frac{1}{2}\right)=1.
\end{equation}
The previous result is summarized in the following 
\begin{propn}\label{prop1}
For $\mu>-\frac{1}{2}$ the Lommel function $s_{\mu,\nu}(z)$ possess the following integral representation:
\begin{equation}\label{intr}
s_{\mu,\nu}(z)=z^{\mu}\int_0^1 \sin(zt)f_{\mu,\nu}(t)dt=z^{\mu}\int_0^1 \frac{\sin(zt)}{(1-t)^{\frac{1}{2}-\mu}}\frac{_2F_1\left(\frac{1}{2}+\nu,\frac{1}{2}-\nu;\mu+\frac{1}{2};\frac{1-t}{2}\right)}{_2F_1\left(\frac{1}{2}+\nu,\frac{1}{2}-\nu;\mu+\frac{1}{2};\frac{1}{2}\right)}dt.
\end{equation}
\end{propn}
Notice that the function $_2F_1\left(\frac{1}{2}+\nu,\frac{1}{2}-\nu;\mu+\frac{1}{2};\frac{1-t}{2}\right)$ is proportional to the Ferrers functions $P_{\nu-\frac{1}{2}}^{-\mu+\frac{1}{2}}(t)$. Indeed, by definition, one has \cite{Pru} 
\begin{equation}
P_{\nu}^{\mu}(t)=\frac{1}{\Gamma(1-\mu)}\left(\frac{1+t}{1-t}\right)^{\frac{\mu}{2}}\, _2F_1\left(-\nu, \nu+1, 1-\mu,\frac{1-t}{2}\right),
\end{equation}
giving
\begin{equation}\label{news}
s_{\mu,\nu}(z)=z^{\mu}\int_0^1 \frac{\sin(zt)}{\left(1-t^2\right)^{\frac{1}{4}-\frac{\mu}{2}}}\frac{P_{\nu-\frac{1}{2}}^{-\mu+\frac{1}{2}}(t)}{P_{\nu-\frac{1}{2}}^{-\mu+\frac{1}{2}}(0)}dt.
\end{equation}
The previous equation can be derived also from a formula given in \cite{Pru}. Indeed one has (see \cite{Pru} equation 2.17.9.1) 
\begin{equation}\label{prud}\begin{split}
&\int_0^1\frac{t^{\alpha-1}}{(1-t^2)^{\frac{\mu}{2}}}\sin(zt)P_{\nu}^{\mu}(t)dt=z\frac{\sqrt{\pi}}{2^{\alpha-\mu+1}}\frac{\Gamma(\alpha+1)}{\Gamma(\frac{2+\alpha-\mu-\nu}{2})\Gamma(\frac{3+\alpha-\mu+\nu}{2})}\, \\
&_2F_3\left(\frac{\alpha+1}{2},\frac{\alpha}{2}+1;\frac{3}{2},\frac{2+\alpha-\mu-\nu}{2},\frac{3+\alpha-\mu+\nu}{2};-\frac{z^2}{4}\right)
\end{split}
\end{equation}
By setting $\alpha=1$ and changing $\mu \to -\mu+\frac{1}{2}$ and $\nu \to \nu-\frac{1}{2}$, equation (\ref{prud}) gives (\ref{news}) via (\ref{iper12}). In this work, however, since we are going to use an old result of Hurwitz \cite{Hurwitz} on the zeros of hypergeometric functions, we prefer to keep the integral (\ref{intr}) as it is. 

We are interested in the monotonicity of the integrand $f_{\mu,\nu}(t)$ for $t\in (0,1)$ and $\mu>-1/2$. We notice that the derivative of $f_{\mu,\nu}$ with respect to $t$ is proportional again to a hypergeometric function:
\begin{equation}\label{deriv}
\frac{df_{\mu,\nu}}{dt}=\frac{1-2\mu}{2(1-t)^{\frac{3}{2}-\mu}}\frac{_2F_1\left(\frac{1}{2}+\nu,\frac{1}{2}-\nu;\mu-\frac{1}{2};\frac{1-t}{2}\right)}{_2F_1\left(\frac{1}{2}+\nu,\frac{1}{2}-\nu;\mu+\frac{1}{2};\frac{1}{2}\right)}.
\end{equation}
The previous relation can be written in a more compact form as
\begin{equation}\label{deriv1}
\frac{df_{\mu,\nu}}{dt}+a_{\mu,\nu}f_{\mu-1,\nu}=0,
\end{equation}
where we set
\begin{equation}\label{amn}
a_{\mu,\nu}\doteq 2\frac{\Gamma\left(\frac{\mu+1+\nu}{2}\right)\Gamma\left(\frac{\mu+1-\nu}{2}\right)}{\Gamma\left(\frac{\mu+\nu}{2}\right.)\Gamma\left(\frac{\mu-\nu}{2}\right)}.
\end{equation}
Notice that $a_{\mu,\nu}$ is related to the integral of $f_{\mu,\nu}(t)$ between $0$ and $1$: indeed, from (\ref{deriv1}) we get for $\mu+1/2>0$
\begin{equation}
\int_0^1f_{\mu,\nu}(t)dt=\frac{1}{a_{\mu+1,\nu}}=\frac{\Gamma\left(\frac{\mu+1+\nu}{2}\right)\Gamma\left(\frac{\mu+1-\nu}{2}\right)}{2\Gamma\left(\frac{\mu+\nu}{2}+1\right.)\Gamma\left(\frac{\mu-\nu}{2}+1\right)}.
\end{equation}
The previous quantity is finite under the given assumptions that $\mu\pm\nu$ is not an odd negative integer.

As regards the monotonicity, from (\ref{deriv1}) we see that the zeros of $\frac{df_{\mu,\nu}}{dt}$ between $0$ and $1/2$ are the zeros of $_2F_1\left(\frac{1}{2}+\nu,\frac{1}{2}-\nu;\mu-\frac{1}{2};\frac{1-t}{2}\right)$, since the constant coefficient is zero only for $\mu+\frac{1}{2}=-n$, $n=0,1,2,...$ as can be seen directly from the explicit expression
\begin{equation}
_2F_1\left(\frac{1}{2}+\nu,\frac{1}{2}-\nu;\mu+\frac{1}{2};\frac{1}{2}\right)=2^{1/2-\mu} \sqrt{\pi} \frac{\Gamma(\mu+\frac{1}{2})}{\Gamma(\frac{\mu+\nu+1}{2})\Gamma(\frac{\mu-\nu+1}{2})} 
\end{equation}
The number of zeros of the hypergeometric functions $ _2F_1(a,b;c;x)$ for $x \in (0,1)$ have been analyzed by  Klein \cite{Klein} and Hurwitz \cite{Hurwitz}.  For completeness, we report these results as given by Hurwitz. They are summarized in Table (\ref{T1}) for $a>b$, $c\leq a+b$\footnote{These choices are not restrictive, since $a$ and $b$ can be interchanged and, for $c>a+b$ one can use the equivalence $_2F_1(a,b;c;x)=(1-x)^{c-a-b}_2F_1(c-a,c-b;c;x)$.}

\begin{figure}
\centering
\includegraphics[scale=0.6]{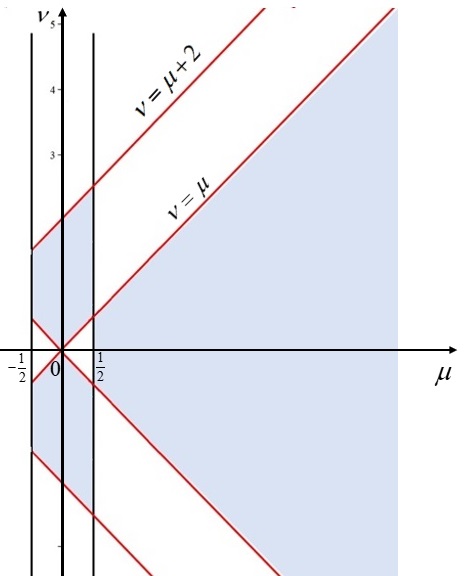}
\caption{The region of monotonicity of $f_{\mu,\nu}(t)$ for $t\in (0,1)$.}
\label{fig1}
\end{figure}

\begin{table}[h]
\begin{center}
\begin{tabular}{|c|c|c|c|}
\hline
$a$ & $b$ & $c$ & $N$\\
\hline
+ & + & + & 0\\
\hline
+ & + & - & $\frac{1+(-1)^{\lfloor -c \rfloor}}{2}$\\
\hline
+ & - & + & $1+\lfloor -b \rfloor$\\
\hline
+& - & - & $\lfloor -b \rfloor - \lfloor -c \rfloor$, for $\lfloor -b \rfloor > \lfloor -c \rfloor$\\
\hline
+& - & - & $\frac{1-(-1)^{\lfloor -b \rfloor + \lfloor -c \rfloor}}{2}$, for $\lfloor -b \rfloor \leq \lfloor -c \rfloor$\\
\hline
-&-&-& $\frac{1+(-1)^{\lfloor -a \rfloor + \lfloor -b \rfloor + \lfloor -c \rfloor}}{2}$\\
\hline
\end{tabular}
\caption{The number of zeros $N$ of $_2F_1(a,b;c;x)$ for $x \in (0,1)$, $a>b$, $c\leq a+b$.}
\label{T1}
\end{center}
\end{table}
We notice however that from table (\ref{T1}) we can provide only partial results about the region of the parameters $(\nu,\mu)$ where the function\\ $_2F_1\left(\frac{1}{2}+\nu,\frac{1}{2}-\nu;\mu-\frac{1}{2};\frac{1-t}{2}\right)$ is free from zeros for $t \in (0,1)$. Indeed, the table  (\ref{T1}) gives this information for $t \in (-1,1)$. We need to extend the results of Hurwitz and Klein and look at the zeros of $ _2F_1(a,b;c;x)$ in the region $x \in (0,1/2)$, i.e. the zeros of $_2F_1\left(\frac{1}{2}+\nu,\frac{1}{2}-\nu;\mu-\frac{1}{2};\frac{1-t}{2}\right)$ for $t \in (0,1)$. To this aim, the quadratic (in the independent variable) functional identities for the hypergeometric functions are very useful. The identity we need is the following (see \cite{Gou}, formula 41 at page 120):
\begin{equation}\label{funct}
\frac{_2F_1\left(\frac{1}{2}+\nu,\frac{1}{2}-\nu;\mu-\frac{1}{2};x\right)}{(1-2x)\left(1-x\right)^{\mu-3/2}}=\, _2F_1\left(\frac{\mu+\nu}{2},\frac{\mu-\nu}{2};\mu-\frac{1}{2};4x(1-x)\right)
\end{equation}
where $x \in (0,1/2)$.  Equation (\ref{funct}) provides an explicit relation between the zeros of $_2F_1\left(\frac{\mu+\nu}{2},\frac{\mu-\nu}{2};\mu-\frac{1}{2};x\right)$ and the zeros of $_2F_1\left(\frac{1}{2}+\nu,\frac{1}{2}-\nu;\mu-\frac{1}{2};x\right)$. In particular  it follows that if $_2F_1\left(\frac{\mu+\nu}{2},\frac{\mu-\nu}{2};\mu-\frac{1}{2};x\right)$ is free from zeros for $x\in(0,1)$, $_2F_1\left(\frac{1}{2}+\nu,\frac{1}{2}-\nu;\mu-\frac{1}{2};x\right)$ is free from zeros for $x\in (0,1/2)$. 
By looking at table (\ref{T1}) for $_2F_1\left(\frac{\mu+\nu}{2},\frac{\mu-\nu}{2};\mu-\frac{1}{2};x\right)$ we get the following
\begin{propn}\label{prop1}
For $\mu >-1/2$ the function $f_{\mu,\nu}(t)$ is monotonic for $t\in (0,1)$ iff the following conditions on the parameters are satisfied:
\begin{itemize}
\item $\mu >1/2$ and $|\nu| \in (0,\mu)$. 
\item $\mu <1/2$ and $ |\nu| \in (|\mu|, \mu+2)$
\end{itemize}
\end{propn}
The region in the plane $(\mu,\nu)$ corresponding to monotonic function $f_{\mu,\nu}(t)$ for $t\in (0,1)$  is illustrated in figure (\ref{fig1}).
The function $f_{\mu,\nu}(t)$ can be decreasing or increasing for $t\in (0,1)$. Since $f_{\mu,\nu}(0)=1$, by looking at what happens in a neighborhood of $t=1$ it is simple to get the following
\begin{propn}\label{prop2}
For $\mu >-1/2$ the function $f_{\mu,\nu}(t)$ has the following monotonicity properties going from $t=0$ to $t=1$: 
\begin{itemize}
\item For $\mu >1/2$ and $|\nu| \in (0,\mu)$ the function is decreasing from 1 to 0.
\item For $\mu <1/2$ and $|\nu| \in (|\mu|, \mu+1)$ the function is increasing from 1 to $+\infty$.
\item For $\mu <1/2$ and $|\nu| \in (\mu+1, \mu+2)$ the function is decreasing from 1 to $-\infty$.
\end{itemize}
\end{propn}
From Proposition (\ref{prop2}) and Theorem (\ref{ThP}) we immediately get the following
\begin{cor}\label{cor1}
Apart the branch point at $z=0$, the function $s_{\mu,\nu}(z)$ for $\mu \in (-1/2,1/2)$ and  $|\nu| \in (|\mu|, \mu+1)$ possesses only real zeros. The zeros are simple and the intervals $(k\pi, (k+1)\pi)$, $k=1, 2, \ldots$ contain the non-negative zeros, each interval containing just one zero.
\end{cor}
The previous generalizes the Corollary (\ref{cor0}) given in \cite{Z}. Further, since for $\mu >1/2$ and  $|\nu| \in (0,\mu)$ the function $f_{\mu,\nu}(t)$ is positive and decreasing, we get the Corollary:
\begin{cor}\label{cor2}
The function $s_{\mu,\nu}(z)$ for $\mu >1/2$ and $|\nu| \in (0,\mu)$ is positive on the positive real axis and possesses only complex zeros.
\end{cor}
Corollary (\ref{cor2}) is not new, it has been given by \cite{Steinig} and then appears in \cite{Gasper}. The proof given here is however simpler. 
The zeros of $s_{\mu,\nu}(z)$ can be better characterized by looking at the asymptotic expansion of the hypergeometric $_1F_2$ function, related to $s_{\mu,\nu}(z)$ by the formula (\ref{iper12}). Indeed, by using the results given in \cite{Lin} and with the aid of (\ref{iper12}) we get, for $z$ real and positive, the asymptotic relation
\begin{equation}\label{asrel}
\frac{s_{\mu,\nu}(z)}{z^{\mu+1}} \sim \frac{1}{z^2}\left(1+O(\frac{1}{z})\right)+\frac{\Gamma\left(\frac{\mu+\nu+1}{2}\right)\Gamma\left(\frac{\mu-\nu+1}{2}\right)}{4\sqrt{\pi}}\left(\frac{2}{z}\right)^{\mu+\frac{3}{2}}\cos\left(z-\frac{\pi}{2}(\mu+\frac{3}{2})\right)(1+O(\frac{1}{z}))
\end{equation}
and if $\mu<1/2$ we see that the zeros are asymptotically given by $z_k \sim \pi\left(k+\frac{2\mu+5}{4}\right)$, for $k$ suitable large integers.

The behavior of $f_{\mu,\nu}(t)$ for $t\in (0,1)$, $\mu \in (-1/2, 1/2)$ and $|\nu| \in (\mu+1, \mu+2)$ gives the following result: yhe function $c-f_{\mu,\nu}(t)$ is positive and increasing for $t \in (0,1)$ for any $c\geq 1$. Indeed, if we introduce the function
\begin{equation}
V_{\mu,\nu}(z) \doteq z^\mu \int_{0}^1(c-f_{\mu,\nu}(t))\sin(zt)dt=z^{\mu-1}c(1-\cos(z))-s_{\mu,\nu}(z),
\end{equation}
we can apply the P\'olya Theorem (\ref{ThP}) and get the following
\begin{cor}\label{corv1}
For any $c \geq 1$, the function $V_{\mu,\nu}(z)=c z^{\mu-1}(1-\cos(z))-s_{\mu,\nu}(z)$ for  $\mu \in (-1/2, 1/2)$ and $|\nu| \in (\mu+1, \mu+2)$ possesses only real zeros. The zeros are simple and the intervals $(k\pi, (k+1)\pi)$, $k=1, 2, \ldots$ contain the non-negative zeros, each interval containing just one zero.
\end{cor}

Actually, (\ref{intr}) is not the only integral representation in terms of a trigonometric kernel. Indeed, it is also possible to give a cosine integral representation. Let us assume $\mu>1/2$. Then, by integrating by parts equation (\ref{intr}) we get
\begin{equation}
s_{\mu,\nu}=z^{\mu}\left(\frac{1}{z}-\frac{a_{\mu,\nu}}{z}\int_0^1 \cos(zt)f_{\mu-1,\nu}(t)\right),
\end{equation}
where we used equation (\ref{deriv1}). Let us set
\begin{equation}\label{cmn}
c_{\mu,\nu}\doteq z^{\mu+1}a_{\mu+1,\nu}\int_0^1 \cos(zt)f_{\mu,\nu}(t).
\end{equation}
Equation (\ref{cmn}) is well defined for $\mu>-1/2$. Again, by integrating by parts (\ref{cmn}) and by using the equation $a_{\mu+1,\nu}a_{\mu,\nu}=(\mu^2-\nu^2)$ we get, for $\mu>+1/2$
\begin{equation}\label{c1}
c_{\mu,\nu}= z^{\mu}\left(\mu^2-\nu^2\right)\int_0^1 \sin(zt)f_{\mu-1,\nu}(t)=(\mu^2-\nu^2)zs_{\mu-1,\nu}.
\end{equation}
From (\ref{c1}) it follows that, for $\mu>-1/2$
\begin{equation}\label{smnnn}
s_{\mu,\nu}=\frac{1}{z((\mu+1)^2-\nu^2)}c_{\mu+1,\nu}= \frac{a_{\mu+2,\nu}}{((\mu+1)^2-\nu^2)}z^{\mu+1}\int_0^1 \cos(zt)f_{\mu+1,\nu}(t).
\end{equation}
Notice that the previous integral representation is still convergent for $\mu>-3/2$. It is also possible to check directly, by using the series for the hypergeometric function around $t=1$ that indeed the integral in (\ref{smnnn}) gives $s_{\mu,\nu}$ also for $\mu \in (-3/2,-1/2)$, so we get the following  
\begin{propn}\label{prop1}
For $\mu>-\frac{3}{2}$ the Lommel function $s_{\mu,\nu}(z)$ possesses the following integral representation:
\begin{equation}\label{intrcos}
s_{\mu,\nu}(z)= \frac{a_{\mu+2,\nu}}{((\mu+1)^2-\nu^2)}z^{\mu+1}\int_0^1 \cos(zt)f_{\mu+1,\nu}(t)=\frac{z^{\mu+1}}{a_{\mu+1,\nu}}\int_0^1 \cos(zt)f_{\mu+1,\nu}(t).
\end{equation}
where $a_{\mu,\nu}$ is given in equation (\ref{amn}) and $f_{\mu,\nu}(t)$ by equation (\ref{intr}).
\end{propn}
It is possible to get the integral representation (\ref{intrcos}) also by directly making the ansatz $s_{\mu,\nu}=z^{\mu+1}\int_0^1\cos(zt)f(t)$ for some $f(t)$ and then by looking at the differential equation and boundary conditions that $f(t)$ must obey so that $s_{\mu,\nu}$ satisfies equation (\ref{Lom}), like we did with the integral representation (\ref{intr}). Since in the integral (\ref{intrcos}) it appears the function $f_{\mu+1,\nu}(t)$, i.e. the same function appearing in the integral (\ref{intr}) but with $\mu \to \mu+1$, from Proposition (\ref{prop1}) we get directly the following 
\begin{propn}
For $\mu >-3/2$ the function $f_{\mu+1,\nu}(t)$ is monotonic for $t\in (0,1)$ iff the following conditions on the parameters are satisfied:
\begin{itemize}
\item $\mu >-1/2$ and  $|\nu| \in (0,\mu+1)$. 
\item $\mu <-1/2$ and $ |\nu| \in (|\mu+1|, \mu+3)$
\end{itemize}
\end{propn}
To understand the behavior of $\frac{f_{\mu+1,\nu}(t)}{a_{\mu+1,\nu}}$, we look at the end-points $t=0$ and $t=1$. For $t=0$ we get
\begin{equation}\label{f0}
 \frac{f_{\mu+1,\nu}(0)}{a_{\mu+1,\nu}}=\frac{1}{a_{\mu+1,\nu}}=\frac{\Gamma\left(\frac{\mu+1+\nu}{2}\right)\Gamma\left(\frac{\mu+1-\nu}{2}\right)}{2\Gamma\left(1+\frac{\mu+\nu}{2}\right.)\Gamma\left(1+\frac{\mu-\nu}{2}\right)}
\end{equation}
For $\mu>-1/2$ $f_{\mu+1,\nu}(1)=0$, giving a decreasing function for $ |\nu| \in (0,\mu+1)$ since the values of the function (\ref{f0}) is positive in this interval. For $\mu<-1/2$ $f_{\mu+1,\nu}(t)$ diverges at $t=1$, the sign being negative. Also, for $|\nu| \in (|\mu+1|, \mu+2)$ the values of the function (\ref{f0}) are negative, whereas for $|\nu| \in (\mu+2, \mu+3)$ are positive. It follows that the function
\begin{equation}
-\frac{1}{a_{\mu+1,\nu}}f_{\mu+1,\nu}(t)
\end{equation}
is positive and increasing for $\mu \in (-3/2, -1/2)$ and $|\nu| \in (|\mu+1|, \mu+2)$. The previous results let to expand the region in the $(\mu,\nu)$ plane for the Corollary (\ref{cor1}). Indeed we immediately get the following 
\begin{cor}\label{cor4}
Apart the branch point at $z=0$, the function $s_{\mu,\nu}(z)$ for $\mu \in (-3/2,-1/2)$ and $|\nu| \in (|\mu+1|, \mu+2)$ possesses only real zeros. The zeros are simple and the intervals $(\frac{(2k+1)\pi}{2}, \frac{(2k+3)\pi}{2})$, $k=0, 1, 2, \ldots$ contain the non-negative zeros, each interval containing just one zero.
\end{cor}
With the same considerations given after Corollary (\ref{cor2}) it is possible to show that the zeros in Corollary (\ref{cor4}) are asymptotically given by $z_k \sim \pi\left(k+\frac{2\mu+5}{4}\right)$ for large $k$.
\begin{figure}
\centering
\includegraphics[scale=0.8]{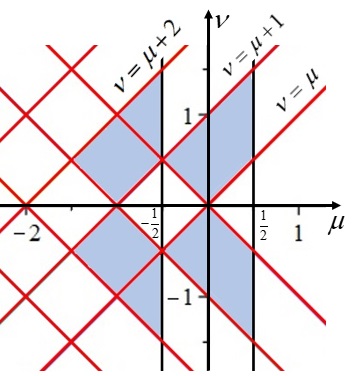}
\caption{The regions described in Corollary (\ref{cor1}) and Corollary (\ref{cor4}).}
\label{fig2}
\end{figure}
As we did with Corollary (\ref{corv1}), we can introduce the function 
\begin{equation}\label{Umunu}
U_{\mu,\nu}(z)=z^{\mu}(c+\frac{1}{a_{\mu+1,\nu}})\sin(z)-s_{\mu,\nu}(z),
\end{equation}
and, from the P\'olya Theorem (\ref{ThP}) we get the following
\begin{cor}\label{corv2}
For any $c \geq 0$, the function $U_{\mu,\nu}(z)$ (\ref{Umunu}) for  $\mu \in (-3/2, -1/2)$ and $|\nu|\in (\mu+2, \mu+3)$ or $\mu>-1/2$ and $|\nu|\in (0,\mu+1)$ possesses only real zeros. The zeros are simple and the intervals $(\frac{(2k+1)\pi}{2}, \frac{(2k+3)\pi}{2})$, $k=0, 1, 2, \ldots$ contain the non-negative zeros, each interval containing just one zero.
\end{cor}
The Corollaries (\ref{corv1}) and (\ref{corv2}) are trivial when the value of $c$ is large and for large real values of $z$, due to the asymptotic (\ref{asrel}). They are less obvious for smaller values of $z$.

By looking closer at the paper \cite{P}, one sees that Theorem (\ref{ThP}) can be actually extended to linear combinations of the functions $U(z)$ and $V(z)$ as also underline by P\'olya himself (see also \cite{Szego}). Under the same assumptions about $f(t)$ of Theorem (\ref{ThP}) indeed, one has \cite{Szego} that the entire function
\begin{equation}
\cos(\theta)U(z)+\sin(\theta)V(z)=\int_0^1f(t)\cos(zt-\theta)dt, \quad \theta \in (0,\pi)
\end{equation}
possesses only real simple zeros, each zero belonging to the intervals
\begin{equation}\label{genint}
\left(\left(k-\frac{1}{2}\right)\pi+\theta,\left(k+\frac{1}{2}\right)\pi+\theta\right),\quad k=0,\pm 1, \pm 2, \ldots
\end{equation}
Due to the results given in Proposition (\ref{prop2}), the Corollaries  (\ref{cor1}) and (\ref{cor4}) can be summarized in the following
\begin{cor}\label{corfin}
The function $z^{-\mu}\left(a_{\mu,\nu}\cos(\theta)s_{\mu-1,\nu}(z)+\sin(\theta)s_{\mu,\nu}(z)\right)$ for $\mu \in (-1/2,1/2)$ and  $|\nu| \in (|\mu|, \mu+1)$ possesses only real zeros for any $\theta \in [0,\pi]$. The zeros are simple and each of the intervals (\ref{genint}) contain just one zero.
\end{cor}

Let us finally notice that for particular values of the parameters $\mu$ and $\nu$ it is possible to get algebraic functions for the kernel of the integral (\ref{intr}). Also, when $\mu$ is a positive integer, it is possible to give explicit formulae for the hypergeometric kernel in terms of trigonometric functions. These observations will be further investigated in a separate paper.

\section{Some properties of the hypergeometric $_1F_2$ function}
P\'olya \cite{P} and Hille \cite{Hille} have been two of the few authors to investigate  about the distribution of the zeros of the hypergeometric function $_1F_2(1;b,c;z)$. For uniformity of notation, we use the same set of parameters for $_1F_2$ as in (\ref{iper12}). The results of P\'olya and Hille can be summarized as follows:
\begin{propn} 
The function $_1F_2\left(1;\frac{\mu-\nu+3}{2},\frac{\mu+\nu+3}{2};z\right)$ possesses only real zeros for $\nu=1/2$ and $\mu \in (-5/2,-1/2)$ or $\mu \in (-7/2,-5/2)$. For $\nu=1/2$ and $\mu>1/2$ it has only complex zeros. For $\nu=1/2$ and $\mu<1/2$ there are infinitely many real zeros.
\end{propn}
Very recently Sokal \cite{Sokal} gives certain conditions ensuring that $_pF_q (a_1,...a_p;b_1,...b_q;x)$, $p\leq q$ possesses only real, non positive roots. In particular he showed that $_pF_q (a_1,...a_p;b_1,...b_q;x)$, $p\leq q$ are entire functions belonging to the Laguerre–P\'olya class $LP^+$ (i.e. they can be obtained as a limit, uniformly on compact subsets of $\mathbb{C}$, of a sequence of polynomials with roots on $(-\infty, 0]$) for arbitrarily large $b_{p+1},...,b_{q}$ if and only if, after a possible reordering, the differences $a_i -b_i$, $i=1,...,p$, are nonnegative integers. When the differences $a_i-b_i$, $i=1,...,p$, are not integers it is still possible for the functions $_pF_q$ to belong to the Laguerre-P\'olya class $LP^+$ for some finite values of $b_{p+1},...b_q$, as it is shown also by the results of P\'olya and Hille stated above. 

As far as we know, if one considers the function $_1F_2(a_1;b_1,b_2;z)$, there are no results in the literature about two or three dimensional sets of $(a_1,b_1,b_2) \in \mathbb{R}^{3}$ for which $_1F_2(a_1;b_1,b_2;z)$ belong to $LP^+$ (see also at the end of \cite{Sokal}). However,  due to the equivalence (\ref{iper12}), any result about the zeros of $s_{\mu,\nu}(z)$ can be directly transferred to $_1F_2\left(1;\frac{\mu-\nu+3}{2},\frac{\mu+\nu+3}{2};z\right)$. Actually, it would be preferable to investigate the zeros of the function $_1F_2\left(1;\frac{\mu-\nu+3}{2},\frac{\mu+\nu+3}{2};z\right)$ rather then those of $s_{\mu,\nu}(z)$ since $_1F_2\left(1;\frac{\mu-\nu+3}{2},\frac{\mu+\nu+3}{2};z\right)$ is an entire function of $z$ for any choice of the parameters $(\mu, \nu)$. In this section, for completeness, we report the Corollaries for $_1F_2\left(1;\frac{\mu-\nu+3}{2},\frac{\mu+\nu+3}{2};z\right)$ corresponding to the Corollaries (\ref{cor1}), (\ref{cor2}), (\ref{corv1}), (\ref{cor4}) and (\ref{corv2}) for $s_{\mu,\nu}(z)$. Also, for completeness, we will give the corresponding regions in the plane $(b_1,b_2)$ for which $_1F_2(a_1;b_1,b_2;z)$ belong to $LP^+$.

\begin{cor}\label{corhy1}
The function $_1F_2\left(1;\frac{\mu-\nu+3}{2},\frac{\mu+\nu+3}{2};z\right)$ for $\mu \in (-1/2,1/2)$ and $|\nu| \in (|\mu|, \mu+1)$ possesses only real negative zeros. The zeros are simple and are contained in the intervals $(-\frac{k^2\pi^2}{4}, -\frac{(k+1)^2\pi^2}{4})$, $k=1, 2, \ldots$, each interval containing just one zero.
\end{cor}

\begin{cor}\label{corhy4}
The function $_1F_2\left(1;\frac{\mu-\nu+3}{2},\frac{\mu+\nu+3}{2};z\right)$  for $\mu \in (-3/2,-1/2)$ and $|\nu| \in (|\mu+1|, \mu+2)$ possesses only real negative zeros. The zeros are simple and are contained in the intervals $(-\frac{(2k+1)^2\pi^2}{16}, -\frac{(2k+3)^2\pi^2}{16})$, $k=0, 1, 2, \ldots$, each interval containing just one zero.
\end{cor}

From the previous two Corollaries, by setting $a_1=1$, $b_1=\frac{\mu-nu+3}{2}$, $b_2=\frac{\mu+\nu+3}{2}$ and by noticing that $_1F_2(a_1;b_1,b_2;z)$ is an entire function of fractional order of growth and hence possesses an infinite number of zeros \cite{titchmarsh}, we get the following
\begin{propn}\label{reg}
In the regions enclosed by the lines $b_1=\frac{1}{2}+\delta$, $b_1=1+\delta$, $b_2=1+\delta$ and $b_1+b_2=\frac{5}{2}+\delta$, where $\delta=0, \frac{1}{2}$, and in the regions enclosed by the lines obtained by the exchange $b_1 \leftrightarrow b_2$, the function $_1F_2(a_1;b_1,b_2;z)$ belong to the Laguerre-P\'olya class $LP^+$. 
\end{propn}
The regions described in Proposition (\ref{reg}) are reported in the following figure. Clearly it is a simple rotation and translation of the regions in figure (\ref{fig2}).
\begin{figure}[H]
\centering
\includegraphics[scale=0.6]{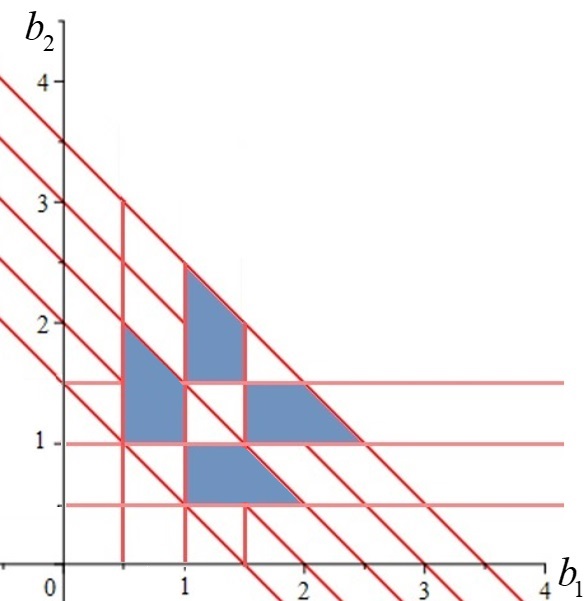}
\caption{The regions (in grey) described in Proposition (\ref{reg}).}
\label{fig3}
\end{figure}

Finally, we give the following results as consequences of the Corollaries (\ref{cor2}), (\ref{corv1}) and (\ref{corv2}).

\begin{cor}\label{corhy2}
The function $_1F_2\left(1;\frac{\mu-\nu+3}{2},\frac{\mu+\nu+3}{2};z\right)$ for $\mu >1/2$ and $|\nu| \in (0,\mu)$ is positive on the positive real axis and possesses only complex zeros.
\end{cor}

\begin{cor}\label{corvhy1}
For any $c \geq 1$, the function 
\begin{equation}
c \left((\mu+1)^2-\nu^2\right)\frac{(1-\cos(z))}{z^2} -\, _1F_2\left(1;\frac{\mu-\nu+3}{2},\frac{\mu+\nu+3}{2};-\frac{z^2}{4}\right)
\end{equation}
 for  $\mu \in (-1/2, 1/2) $ and $|\nu| \in (\mu+1, \mu+2)$ possesses only real zeros. The zeros are simple and the intervals $(k\pi, (k+1)\pi)$, $k=1, 2, \ldots$ contain the non-negative zeros, each interval containing just one zero.
\end{cor}

\begin{cor}\label{corvhy2}
For any $c \geq 0$, the function 
\begin{equation}
a_{\mu+2,\nu}(c a_{\mu+1,\nu}+1)\frac{\sin(z)}{z} - \, _1F_2\left(1;\frac{\mu-\nu+3}{2},\frac{\mu+\nu+3}{2};-\frac{z^2}{4}\right)
\end{equation}
 for  $\mu \in (-3/2, -1/2)$ and $\nu \in (\mu+2, \mu+3)$ or $\mu>-1/2$ and $\nu\in (0,\mu+1)$ possesses only real zeros. The zeros are simple and the intervals $(\frac{(2k+1)\pi}{2}, \frac{(2k+3)\pi}{2})$, $k=0, 1, 2, \ldots$ contain the non-negative zeros, each interval containing just one zero.
\end{cor}

\begin{center} {\bf Acknowledgments} \end{center}
I wish to acknowledge the support of Universit\`a degli Studi di Brescia,  GNFM-INdAM and INFN, Gr. IV - Mathematical Methods in NonLinear Physics.
\bigskip




\end{document}